




\input amstex
\documentstyle{amsppt}
\magnification=\magstep1
\pagewidth{6truein}
\pageheight{8.8truein}
\NoBlackBoxes
\def\A{{\Cal A}}
\def\D{{\Cal D}}
\def\F{{\Cal F}}
\def\nat{{\Bbb N}}
\def\que{{\Bbb Q}}
\def\real{{\Bbb R}}
\def\ep{{\varepsilon}}
\def\To{\Rightarrow}
\def\supp{\operatorname{supp}}
\def\dis{\ {\buildrel {\D}\over =}\ }
\def\normm{{|\!|\!|\,}}
\def\dotnorm{|\cdot|}
\def\dotNorm{\|\cdot\|}
\def\dotnormm{\normm\cdot\normm}
\topmatter
\title A problem on spreading models\endtitle
\author E. Odell and Th. Schlumprecht\endauthor
\abstract It is proved that if a Banach space $X$ has a basis $(e_n)$ 
satisfying every spreading model of a normalized block basis of $(e_n)$ 
is 1-equivalent to the unit vector basis of $\ell_1$ (respectively, $c_0$) 
then $X$ contains $\ell_1$ (respectively, $c_0$).
Furthermore Tsirelson's space $T$ is shown to have the property that 
every infinite dimensional subspace contains a sequence having spreading 
model 1-equivalent to the unit vector basis of $\ell_1$. 
An equivalent norm is constructed on $T$ so that $\|s_1+s_2\|<2$ 
whenever $(s_n)$ is a spreading model of a normalized basic sequence in $T$. 
\endabstract
\thanks Research supported by NSF and TARP.\endthanks 
\endtopmatter

\document  
\baselineskip=18pt			
\subhead \S0. Introduction\endsubhead 

From the fact that $\ell_1$ (and $c_0$) are not distortable \cite{J} it 
follows that if a Banach space $X$ contains $\ell_1$ (or $c_0$) then some 
basic sequence $(e_i)$ in $X$ has the property that every spreading model 
of a normalized block basis is 1-equivalent to the unit vector basis of 
$\ell_1$ (or $c_0$). 
In this paper we prove the converse statements. 
More generally we show 

\proclaim{Theorem A}
Let $(e_i)$ be a basis for $X$ 
\roster
\item"a)" If $\|s_1+s_2\|=2$ whenever $(s_n)$ is any spreading model of a 
normalized block basis of $(e_i)$, then $X$ contains a subspace isomorphic 
to $\ell_1$. 
\item"b)" If $\|s_1+s_2\|=1$ whenever $(s_n)$ is any spreading model of a 
normalized block basis of $(e_i)$, then $X$ contains a subspace isomorphic 
to $c_0$.
\endroster
\endproclaim 

The proof of a) will be achieved by showing that such an $(e_i)$ cannot 
be weakly null, if normalized. 
The proof will depend heavily on the theory of the generalized Schreier 
classes of subsets of $\nat$ as introduced in \cite{AA}. 
We make strong use of recent results in \cite{AMT} as well as a result from 
\cite{AO}. 

From Theorem A we obtain the 

\proclaim{Corollary} 
Let $(e_i)$ be a basis for $X$. 
If $X$ contains no subspace isomorphic to $\ell_1$ or $c_0$ then there 
exists a normalized block basis of $(e_i)$ having spreading model $(s_i)$ 
satisfying 
$$1< \|s_1 + s_2\| <2\ .$$ 
\endproclaim 

Theorem A a) is proved in \S2 while part b) and the (easy) Corollary 
are proved in \S3. 

In conjunction with Theorem A it is worth considering Tsirelson's space $T$. 
$T$ is reflexive with an unconditional basis $(t_i)$ and yet all spreading 
models of normalized block bases of $(t_i)$ are 2-equivalent to the unit 
vector basis of $\ell_1$. 
Furthermore we have 

\proclaim{Theorem B} 
Let $X$ be an infinite dimensional subspace of $T$. 
Then there exists $(x_i)\subseteq X$ with spreading model 1-equivalent to 
the unit vector basis of $\ell_1$. 
\endproclaim 

We prove this theorem in \S4. 
Furthermore we show that $T$ can be renormed to fail the conclusion 
of Theorem~B. 

We do not know if Theorem A can be extended to $\ell_p$ ($1<p<\infty$). 

\definition{Problem}
Let $(e_i)$ be a basis for $X$ and $1<p<\infty$. 
Suppose that every spreading model of any normalized block basis of $(e_i)$ 
is 1-equivalent to the unit vector basis of $\ell_p$. 
Does $X$ contain $\ell_p$, either isomorphically or almost isometrically?
\enddefinition 

\subhead \S1 Preliminaries\endsubhead 

\definition{Definition 1.1} 
Let $(e_i)$ be a normalized basic sequence. 
A basic sequence $(s_i)$ is a {\it spreading model\/} of $(e_i)$ if 
for some sequence $\ep_n\downarrow 0$  and all $(a_i)_1^n \subseteq 
[-1,1]^n$ we have 
$$\left| \Big\| \sum_1^n a_i e_{k_i}\Big\| - \Big\| \sum_1^n a_i s_i
\Big\| \right| < \ep_n\text{ whenever } n\le k_1 <\cdots < k_n\ .$$ 
\enddefinition 

It is well known that every normalized basic sequence has a subsequence 
with a spreading model. 
Also $(s_i)$ is necessarily {\it spreading\/}
($\|\sum_1^n a_i s_i\| = \|\sum_1^n 
a_i s_{k_i}\|$ if $k_1<k_2<\cdots <k_n$). 
If $(e_i)$ is weakly null then $(s_i)$ is suppression-1-unconditional 
($\|\sum_F a_i s_i\| \le \|\sum_1^n a_i s_i\|$ if 
$F\subseteq \{1,\ldots,n\}$). 
These and other results on spreading models can be found in \cite{BL}. 

$[\nat]$ denotes the set of all subsequences of $[\nat]$. 
If $M\in [\nat]$, $[M]$ is the set of all subsequences of $M$. 
$[M]^{<\omega}$ is the class of all finite subsets of $M$. 
If $E,F\in [\nat]^{<\omega}$, ``$E<F$'' means that $\max E <\min F$, 
``$k<E$'' means $\{k\} <E$.

\definition{Definition 1.2} 
\cite{AA} {\bf Generalized Schreier classes}.   
The classes $(S_\alpha)_{\alpha<\omega_1}$ of collections of finite subsets 
of $\nat$ are inductively defined as follows 
$$\align 
S_0 & = \big\{ \{n\} : n\in\nat\big\} \cr 
S_{\alpha+1} & = \Big\{ E : E = \bigcup_1^k E_i \text{ for some } k\in\nat 
\text{ and }\cr
&\qquad k\le E_1 < \cdots < E_k \text{ where } E_i\in S_\alpha \text{ for } 
i\le k\Big\} 
\endalign$$
If $\alpha$ is a limit ordinal, choose $\alpha_n\uparrow \alpha$ and 
set $S_\alpha = \{E:k\le E\in S_{\alpha_k}$ for some $k\in\nat\}$. 
We also consider the empty set  $\emptyset \in S_\alpha$ for all $\alpha$. 
The definition of $S_\alpha$ depends upon this particular choice of 
$(\alpha_n)$ but the results we use concerning the $S_\alpha$'s are 
independent of that choice. 
\enddefinition 

The Schreier classes have played a prominent role in a number of recent 
papers (e.g., \cite{AA}, \cite{AD}, \cite{AMT}, \cite{OTW}, \cite{AO}). 

If $M= (m_i)\in[\nat]$ and $\alpha <\omega_1$, $S_\alpha (M) = \{(m_i)_{i\in F}
: F\in S_\alpha\}$. 
It is easy to see that $S_\alpha (M) \subseteq S_\alpha$. 
We also recall that the classes $S_\alpha$ (or $S_\alpha (M)$) are all 
{\it regular\/}. By this we mean they are pointwise closed, {\it hereditary\/} 
$(E\subseteq F\in S_\alpha \To E\in S_\alpha$) and {\it spreading\/} 
($(n_i)_1^k \in S_\alpha$ and $m_1<\cdots < m_k$ with $m_i \ge n_i$ for 
$i\le k$ implies that $(m_i)_1^k \in S_\alpha$). 

\proclaim{Proposition 1.3 \cite{AO}} 
Let $N\in [\nat]$. There exists $M\in [N]$ such that for all $\alpha<\omega_1$ 
if $F\in S_\alpha$ and $F\subseteq M$ then $F\setminus \min (F) \in S_\alpha 
(N)$. 
\endproclaim 

We also need some definitions and a result from \cite{AMT}. 
Let 

\definition{Definition 1.4 \cite{AMT}} 
If $M= (m_i) \in [\nat]$ and $(e_i)$ is a normalized basic sequence we 
inductively define $\alpha_n^M = \alpha_n^M (e_i)$ for $\alpha<\omega_1$ 
and $n\in\nat$ as follows 
$$\align
0_n^M &= e_{m_n}\cr 
(\alpha+1)_1^M &= \frac1{m_1} \sum_{i=1}^{m_1} \alpha_i^M \cr 
(\alpha+1)_{n+1}^M &= \frac1{m_{k_n}} \sum_{i=1}^{m_{k_n}} \alpha_i^{M'}
\endalign$$ 
where $M' = \{m\in M :m>\supp (\alpha+1)_n^M\}$ and $m_{k_n} = \min M'$. 
If $\alpha = \lim \alpha_n$ is a limit ordinal we set 
$$\align 
\alpha_1^M & = (\alpha_{m_1})_1^M\ \text{ and for }\ k>1\ ,\quad 
\alpha_k^M = (\alpha_{n_k})_1^{M_k}\cr 
\text{where }\ M_k & = \{m\in M :m>\supp \alpha_{k-1}^M\}\text{ and }
n_k = \min M_k\ .
\endalign$$
\enddefinition 

\proclaim{Proposition 1.5 \cite{AMT}} 
\roster
\item"1)" For $\alpha <\omega_1$, $(\alpha_n^M) 
= (\alpha_n^M(e_i))_{n=1}^\infty$ is a convex block basis of $(e_i)$ with 
$M= \bigcup_n \supp (\alpha_n^M)$. 
Moreover $\supp \alpha_n^M \in S_\alpha$ for all $n$. 
\item"2)" If $M\in [\nat]$, $\alpha <\omega_1$ and $(n_k)\in[\nat]$ then 
$\alpha_{n_k}^M = \alpha_k^{M'}$ where $M'= \bigcup_k \supp 
(\alpha_{n_k}^M)$. 
\endroster
\endproclaim 

If $x= \sum a_ie_i\in \langle e_i\rangle$ and $F\subseteq \nat$ we 
define $\langle x,F\rangle = \sum_{i\in F} a_i$. 

\definition{Definition 1.6 \cite{AMT}} 
Let $\F$ be an hereditary collection of subsets of $\nat$. 
Let $M\in [\nat]$, $\ep >0$, $\alpha <\omega_1$ 
and let $(e_i)$ be a normalized basic 
sequence. $\F$ is {\it $(M,\alpha,\ep)$ large\/} if for all $N\in [M]$ 
and $n\in\nat$ for $\alpha_n^N = \alpha_n^N (e_i)$,  
$$\sup_{F\in\F} \langle \alpha_n^N, F\rangle >\ep\ .$$ 
\enddefinition 

\proclaim{Theorem 1.7 \cite{AMT, Proposition 2.3.2 and Theorem 2.2.6}} 
If $\F$ is $(M,\alpha,\ep)$ large then there exists $N\in[M]$ with 
$$\F \supseteq S_\alpha (N)\ .$$ 
\endproclaim 

\subhead \S2. $\ell_1$ spreading models\endsubhead 

In this section we prove 

\proclaim{Theorem 2.1} 
Let $(e_i)$ be a normalized basis for $X$ having the property that 
$\|s_1+s_2\|=2$ whenever $(s_n)$ is a spreading model of $(x_i)$ where 
$x_i = y_i/\|y_i\|$ and $(y_i)$ is any convex block subsequence of $(e_i)$ 
satisfying $\lim_i \|y_i\| >0$. 
Then $(e_i)$ is not weakly null. 
Moreover for all $\ep >0$ there exists $M= (m_i)\in [\nat]$ and $x^* \in S 
(X^*)$ with $x^* (e_{m_i}) >1-\ep$ for all $i$. 
\endproclaim 

From Theorem 2.1 it follows that in a space $X$ whose block bases have 
only $\ell_1$ as spreading model no block basis is weakly null. 
Thus 
in light of Rosenthal's theorem \cite{R}, Theorem~A~a) is a quick consequence 
of Theorem~2.1.  

The hypothesis yields that for all $n$, $\|\sum_1^{2^n} s_i\| =2^n$ from 
which it follows that $\|\sum_1^k s_i\|=k$ for all $k$. 
We shall use this below in the following way. 
Given $\ep >0$ there exists a subsequence $(x_{n_k})$ of $(x_i)$ so that 
for all $k$, $\frac1{r_k} \|x_{n_k} + x_{n_{k+1}} +\cdots + x_{n_{k+r_k-1}}\| >
1-\ep$ where $r_k = \min (\supp (x_{n_k}))$, w.r.t. $(e_i)$.  

\demo{Proof} 
Given $\ep>0$ set 
$$\F_\ep = \{F\subseteq \nat : \text{ there exists } x^* \in S(X^*) 
\text{ with } x^* (e_i) >1-\ep \text{ for } i\in F\}\ .$$ 
We shall prove by induction on $\alpha$ that $(P_\alpha)$ holds for all 
$\alpha <\omega_1$ where 
$$(P_\alpha) \qquad \left\{ 
\eqalign{&\text{For all $M\in [\nat]$ and $\ep>0$ there exists}\cr 
&N\in [M]\text{ with } \F_\ep \supseteq S_\alpha (N)\ .\cr}
\right.$$ 

$(P_0)$ is clear. 
If $\alpha$ is a limit ordinal and $S_\alpha$ is defined via the sequence 
$\alpha_n\uparrow \alpha$ we proceed as follows. 
Given $M$ and $\ep>0$ we can choose, by the inductive hypothesis, 
$M\supseteq N_1\supseteq N_2 \supseteq \cdots$ so  that $\F_\ep \supseteq 
S_{\alpha_n} (N_n)$ for $n\in\nat$. 
Let $N_n = (k_i^n)_{i=1}^\infty$ and set $N= (k_n^n)_{n=1}^\infty$. 
Then $\F_\ep \supseteq S_\alpha (N)$. 

Finally assume that $(P_\beta)$ holds and let $M\in [\nat]$ and 
$\alpha = \beta+1$. Let $\ep>0$ and choose $\ep'>0$ so that $\ep' <\ep/2$. 
We may assume that $\F_{\ep'} \supseteq S_\beta (M)$. 
\enddemo

\demo{Claim} 
There exists $N\in [M]$ so that $\F_{3\ep}$ is $(N,\alpha,\ep)$ large.
\enddemo

Indeed for $n\in\nat$ define 
$$\A_n = \Big\{ L\in [M] : \sup_{F\in \F_{3\ep}} \langle \alpha_n^L,F\rangle  
> \ep\Big\}\ .$$ 
Then $\A_n$ is a pointwise closed subset of $[\nat]$ and so $\A = \bigcap_n 
\A_n$ is Ramsey (see \cite{E}, also \cite{O}). 
Thus there exists $N\in [M]$ with either $[N]\subseteq \A$ or $[N]\subseteq 
[\nat]\setminus \A$. 
By passing to a subsequence of $N$ we may assume by Proposition~1.3, that 
for $\gamma <\omega_1$ if $F\in S_\gamma$, $F\subseteq N$ then 
$F\setminus \{\min (F) \in S_\gamma (M)\}$. 
Furthermore we may assume that for all $n$, $\|\beta_n^L \|_\infty<\ep'$ 
for all $L\in [N]$. 
(Indeed this holds if $n_1 = \min N$ satisfies $n_1^{-1} <\ep'$.) 
It follows that for all $n$, since $\supp \beta_n^N \setminus \min (\supp 
\beta_n^N) \in S_\beta (M)\subseteq \F_{\ep'}$, that $\|\beta_n^N\| >
1-2\ep'$. 
From the hypothesis of our theorem applied to $x_n = \beta_n^N/\|\beta_n^N\|$ 
we obtain a subsequence $(\beta_{n_k}^N)_{k=1}^\infty$ satisfying for all $k$: 
$$\frac1{r_k} \|\beta_{n_k}^N + \beta_{n_{k+1}}^N +\cdots +  
\beta_{n_{k+r_k-1}}^N \| > 1-3\ep'
\tag $*$ $$
where $r_k = \min (\supp \beta_{n_k}^N)$. 

\noindent
Let $L= \bigcup_k \supp \beta_{n_k}^N$. 
Then by Proposition~1.5, $\beta_k^L = \beta_{n_k}^N$. 
Hence from $(*)$ and the definition of $\alpha_n^L$ we have $\|\alpha_n^L\| 
> 1-3\ep'$ for all $n$. 

Let $n\in\nat$ and $x^* \in S(X^*)$ with $x^* (\alpha_n^L) > 1-3\ep'$. 
Write $\alpha_n^L = \sum_1^p a_i e_i$ and set $F= \{i :x^* (e_i) >1-3\ep\}$. 
Then $\sum_{i\in F} a_i > \ep$ (otherwise 
$x^* (\alpha_n^L) \le \sum_{i\in F} a_i + 1-3\ep \le 1-2\ep <1-3\ep'$). 
Thus $L\in \A$ and hence $[N]\subseteq \A$, whence the claim follows. 
Thus by Theorem~1.7, $(P_\alpha)$ holds (we actually proved $(P_\alpha)$ 
for $3\ep$ replacing $\ep$). 

Since $(P_\alpha)$ holds for all $\alpha<\omega_1$ we obtain the 
``moreover'' statement of the theorem. 
Indeed this follows easily from an argument of Bourgain \cite{B}. 
Let $T$ be the tree $T= \{ (n_i)_1^k : n_1<\cdots < n_k$ and there exists 
$x^* \in S(X^*)$ with $x^* (e_{n_i}) >1-\ep$ for $i\le k\}$. 
$T$ is a closed tree and thus if $T$ were well founded (no infinite 
branches) then the order of $T$ is $<\omega_1$. 
But since $(P_\alpha)$ holds, the order of $T\ge \omega^\alpha$ for 
all $\alpha$. The latter holds since the order of $S_\alpha (N)$ is 
$\omega^\alpha$, as is well known (see e.g., \cite{AA} or \cite{OTW}).\qed 

\subhead \S3. $c_0$ spreading models\endsubhead 

In this section we prove Theorem A b) and the corollary. 
Note that the hypothesis yields that $\|\sum_1^n s_i\|=1$ for all $n$. 
Also the hypothesis is satisfied if all spreading models of normalized block 
bases of $(e_i)$ are 1-equivalent to the unit vector basis of $c_0$ but this 
is a stronger condition than the hypothesis as the following example 
indicates. 

\example{Example 3.1} 
Let $\|\sum_1^n a_i e_i\| = \max_{i<j} |a_i - a_j|$. 
Then if $X$ is the completion of $(\langle e_i\rangle,\|\cdot\|)$, $X$ 
satisfies the hypothesis of the theorem yet $(e_i)$, which is its own 
spreading model, is not 1-equivalent to the unit vector basis of $c_0$. 
\endexample 

Assume that $X$ has a basis $(e_i)$ satisfying the hypothesis of b). 
We break the proof into several steps. 
For $a,b\in \langle e_i\rangle$ we write ``$a<b$'' if $\supp (a) <\supp (b)$. 

\demo{Step 1} 
For all $\ep>0$ and $\ell\in\nat$ there exists $m\in\nat$ so that for all 
$a>e_m$ with $\|a\|=1$ there exists $b_1<\cdots < b_\ell \le e_m$, 
$\|b_i\| =1$ for $i\le \ell$, such that for all $1\le q\le p\le\ell$ and 
$\delta = 0$ or $1$, 
$$\left| \Big\| \sum_q^p b_i + \delta a\Big\| -1\right| <\ep\ .$$ 
\enddemo 

\demo{Proof} 
If not then there exists $\ep>0$ and $\ell\in\nat$ such that for all 
$m\in\nat$  there exists $a_m>e_m$, $\|a_m\|=1$ so that for all 
$b_1<\cdots < b_\ell \le e_m$, $\|b_i\|=1$ there exists 
$1\le q_m\le p_m\le \ell$ 
with for some $\delta_m = 0$ or $1$, 
$$\left| \Big\| \sum_{q_m}^{p_m} b_i + \delta_m a_m\Big\| -1\right|>\ep\ .$$ 
Choose a subsequence $(a_{m_i})$ of $(a_m)$ having a spreading model $(s_i)$. 
Since $\|\sum_1^{\ell+1} s_i\|=1 = \|\sum_1^\ell s_i\|$, 
we may assume that $\big|\,\| \sum_q^p a_{m_i} +\delta a_{m_{\ell+1}}\| -1
\big| <\ep$ for all $q\le p\le \ell$ and $\delta =0$ or $1$ and that 
$a_{m_\ell} <e_{m_{\ell+1}}$. 
This contradicts our choice of $a_{m_{\ell+1}}$. 

Let $\ep_i \downarrow 0$ with $\sum_1^\infty \ep_i <1$. 
Applying Step~1 to $\ep=\ep_1$ and $\ell=2$ we obtain  $m_1$ so that for all 
$a>e_{m_1}$, $\|a\|=1$ there exist $x_1<y_1 \le e_{m_1}$, 
$\|x_1\| = \|y_1\| =1$ and  
$$\big|\, \|x_1 + y_1 +a\| -1\big| <\ep_1\quad\text{and}\quad 
\big|\,\|y_1 +a\| -1\big| <\ep_1\ .$$ 
\enddemo 

\demo{Step 2} 
There exist $1=m_0 <m_1 <m_2<\cdots$ such that for all $k$ and $a>e_{m_k}$, 
$\|a\| =1$, there exists $x_1<y_1 \le e_{m_1} <x_2 <y_2 \le e_{m_2} 
<\cdots < x_k < y_k \le e_{m_k}$ satisfying 
$\big|\, \|x_i\| -1\big| <\ep_i$, 
$\big|\, \|y_i\|-1\big| <\ep_i$ for $i\le k$ and for all 
$F\subseteq \{1,2,\ldots,k\}$ 
$$\Big\| \sum_F x_i + \sum_1^k y_i +a\Big\| \le 1+\sum_1^k \ep_i\ .$$
\enddemo 

\demo{Proof} 
We proceed by induction on $k$. 
The case $k=1$ was presented above. 
Assume $m_1<\cdots < m_{k-1}$ have been chosen. 
Let $\ep>0$ satisfy ${\ep\over 1-\ep} + (2k-2)\ep<\ep_k$. 
We can, by a compactness argument, find $\ell$ so that if the 
induction hypothesis for $k-1$ is applied to each of $\ell$ different 
$a$'s~$>e_{m_{k-1}}$, say $(a_n)_1^\ell$, then if $(x_i^n)_1^{k-1}, 
(y_i^n)_1^{k-1}$ satisfy Step~2 for $a_n$, for some $n\ne m\le \ell$ we 
have $\|x_i^n - x_i^m\|, \|y_i^n - y_i^m\| <\ep$ for $i\le k-1$. 
Choose $m_k$ by Step~1 applied for this $\ell$ and $\ep$ to 
$(e_i)_{i=m_{k-1}+1}^\infty$.  
Let $a>e_{m_k}$ with $\|a\|=1$. 
Choose $e_{m_{k-1}} <b_1<\cdots < b_\ell \le e_{m_k} \le a$ to satisfy 
$$\left| \Big\| \sum_q^p b_i + \delta a\Big\| -1\right| <\ep\text{ for } 
1\le q\le p \le \ell\ ,\quad \delta = 0,1\ .$$ 
By Step 2 (for $k-1$) there exists for $1\le q\le \ell$, 
$x_1^q < y_1^q \le e_{m_1} <\cdots < x_{k-1}^q <y_{k-1}^q \le 
e_{m_{k-1}}$ so that 
$\big|\,\|x_i^q\|-1\big|,\big|\,\|y_i^q\|-1\big| <\ep_i$ 
for $i\le k-1$ and so that 
for $F\subseteq \{1,\ldots,k-1\}$, 
$$\Big\| \sum_F x_i^q + \sum_1^{k-1} y_i^q + \Big\| \sum_q^\ell b_i+a
\Big\|^{-1} \biggl( \sum_q^\ell b_i +a\biggr)\Big\| 
\le 1+ \sum_1^{k-1} \ep_i\ .$$ 
Thus 
$$\Big\| \sum_F x_i^q + \sum_1^{k-1} y_i^q + \sum_q^\ell b_i +a\Big\| 
< 1+\sum_1^{k-1} \ep_i + {\ep\over 1-\ep}\ .$$
Choose $q'<q$ so that $\|x_i^q - x_i^{q'}\|, \|y_i^q - y_i^{q'}\|<\ep$ for 
$i\le k-1$. 
Let $x_k = \sum_{q'}^{q-1} b_i$ and $y_k = \sum_q^\ell b_i$. 
We have $\big|\,\|x_k\| -1\big|, \big|\, \|y_k\|-1\big| <\ep$ and for 
$F\subseteq  \{1,\ldots,k-1\}$ 
\roster
\item"1)" \quad $\displaystyle \Big\|\sum_F x_i^q + \sum_1^{k-1} 
y_i^q + y_k +a\Big\| < 1+ \sum_1^{k-1} \ep_i 
+ {\ep\over 1-\ep}\ ,$
\item"2)" \quad $\displaystyle \Big\|\sum_F x_i^{q'} + \sum_1^{k-1} 
y_i^{q'} + x_k +y_k +a\Big\| < 1+ \sum_1^{k-1} \ep_i 
+ {\ep\over 1-\ep}\ .$
\endroster 
It follows that if we set $x_i=x_i^q$, $y_i= y_i^q$ for $i\le k-1$ then 
for $F\subseteq \{1,\ldots,k\}$ 
\roster
\item"3)" \quad $\displaystyle \Big\| \sum_F x_i 
+ \sum_1^k y_i +a\Big\| \le 1+\sum_{i=1}^{k-1} \ep_i 
+ {\ep\over 1-\ep} + (2k-2) \ep\ .$ 
\endroster

Thus Step 2 follows by our choice of $\ep$. 

Applying Step 2 to an arbitrary $a_k>e_{m_k}$, $\|a_k\|=1$ we obtain that 
for all $k$ there exists $x_i^k <y_i^k \in \langle e_j
\rangle_{m_{k-1}+1}^{m_k}$ with $\big|\,\| x_i^k\|-1\big| ,\big|\,\|y_i^k\|-1
\big| <\ep_i$ for $i\le k$ and also for $F\subseteq \{1,\ldots,k\}$, 
$$\Big\| \sum_F x_i^k + \sum_1^k y_i^k +a_k\Big\| <2\ .$$ 
It follows that $\|\sum_F x_i^k \|<4$. 

Choose $(k_j)\in [\nat]$ so that for all $i$, $\lim_{j\to\infty} x_i^{k_j}
\equiv x_i$ exists. 
We have for all $F\in [\nat]^{<\omega}$, 
$$\Big\| \sum_F x_i\Big\| \le 4\quad\text{and}\quad 
\big|\, \|x_i\| - 1\big| <1+\ep_i\ .$$ 
Thus $(x_i)$ is equivalent to the unit vector basis of $c_0$.\qed 
\enddemo 

We end this section by presenting the 

\demo{Proof of corollary to Theorem A}

If the corollary is false then all such $s_i$'s satisfy $\|s_1+s_2\| =1$ or 2 
and by Theorem~A both occur. 
Choose $(y_n)$ and $(z_n)$, normalized block bases of $(e_i)$ with spreading 
models $(s_i)$ and $(t_i)$, respectively, satisfying $\|s_1+s_2\|=1$ 
and $\|t_1+t_2\| =2$. 
We may assume that $(y_1,z_1,y_2,z_2,\ldots)$ is a block basis of $(e_i)$. 
Furthermore, by a diagonal argument, we may assume that $(\alpha y_n+\beta 
z_n)_{n=1}^\infty$ has a spreading model $(s_n^{\alpha,\beta})_{n=1}^\infty$ 
for all $\alpha,\beta \in\real$ (not both $0$). 
Now $s_n^{1,0} = s_n$ and $s_n^{0,1} =t_n$. 
There exists a continuous curve $\gamma :[0,1]\to \real^2$, 
$\gamma (t) = (\alpha (t) ,\beta (t))$ so that $\|s^{\alpha(t),\beta(t)}\|=1$ 
for all $t\in [0,1]$ and $\gamma(0)=(1,0)$, $\gamma (1)=(0,1)$. 
We thus obtain by continuity that for all $r\in (1,2)$ there exists $t$ with 
$$\|s_1^{\alpha(t),\beta(t)} + s_2^{\alpha(t),\beta(t)}\|=r\ .
\eqno\qed$$ 
\enddemo 

\remark{Remark 3.2} 
Let $(e_i)$ be a basis for $X$ and let $I(X) = \{r:$ there exists a normalized 
block basis of $(e_i)$ having spreading model $(s_i)$ with $\|s_1+s_2\| =r\}$. 
The proof shows that $I(X)$ is a subinterval of $[1,2]$. 
\endremark

In the next section, Proposition 4.4, we shall see that $I(X)$ need not 
be closed. 
\newpage

\subhead \S4. Theorem B and spreading models of $T$\endsubhead 

If $(e_i)$ is a basic sequence and $x,y\in \langle e_i\rangle$ we say 
that {\it $x$ equals $y$ in distribution\/} ($x\dis y$) if there exist 
$(n(i)),(m(i)) \in [\nat]$ so that $\sum x(i)e_{n(i)} = \sum y(i) e_{m(i)}$. 
The {\it distance in distribution\/} between $x$ and $y$ is defined as 
$$d(x,y) = \inf \big\{ \|\tilde x-\tilde y\| : \tilde x \dis x \text{ and } 
\tilde  y\dis y\big\}\ .$$ 
For $E\subseteq \nat$ we set $Ex = \sum_{i\in E} x(i) e_i$. 

In order to prove Theorem B we first prove 

\proclaim{Proposition 4.1} 
Let $(e_i)$ be a normalized basic sequence having spreading model $(s_i)$ 
which is $K$-equivalent to the unit vector basis of $\ell_1$. 
Assume that there exists a normalized block basis $(x_i)$ of $(s_i)$ which 
satisfies 
\roster
\item"1)" $\lim_i \|x_i\|_{\ell_1} =K$ and 
\item"2)" $(x_i)$ is Cauchy in distribution
\endroster
(i.e., for all $\ep>0$ there exists $n_0$ with $d(x_n,x_m) <\ep$ if 
$n,m\ge n_0$). 
Then there exists a block basis $(y_i)$ of $(e_i)$ having spreading model 
1-equivalent to the unit vector basis of $\ell_1$. 
\endproclaim 

\demo{Proof} 
Let $x_i = \sum_{j\in A_i} a_j^i s_j$ for some choice of scalars and sets 
of integers $A_1<A_2 <\cdots$. 

Using 1), 2) and the fact that $(s_i)$ is $K$-equivalent to the unit vector 
basis of $\ell_1$ we have the following. 
\roster
\item"3)" For all $\ep>0$ there exists $n_0\in\nat$ so that for all 
$n\ge n_0$ there exists $F_n=F_n(\ep) \subseteq A_n$ 
with $|F_n| \le |A_{n_0}|$,  
$\|F_nx_n\|_{\ell_1} >K-\ep$, $d(F_nx_n,x_{n_0}) <\ep$ 
and $\sum_{i\in A_n\setminus F_n} |a_i^n| <\ep$. 
\endroster

For $n\in\nat$ set $y_n = \sum_{i\in A_n} a_i^n e_i$. 
Let $k\in\nat$ and $(a_i)_1^k \subseteq \real$ with $\sum_1^k|a_i|=1$. 
If $n_0<n_1<\cdots <n_k$ then from 3) we obtain 
$$\Big\| \sum_1^k a_i y_{n_i} \Big\| \ge 
\Big\| \sum_1^k a_i F_{n_i} y_{n_i}\Big\| -\ep\ .$$ 
Also 
$$\Big|\bigcup_1^k \supp (F_{n_i} y_{n_i}) \Big| \le k|A_{n_0}| .$$ 
Thus for all $\ep>0$,  
$$\align 
\liminf \Sb n_1\to\infty \\ n_1<\cdots <n_k\endSb 
\Big\| \sum_1^k a_i y_{n_i} \Big\| 
&\ge \liminf \Sb n_1\to\infty\\ n_1<\cdots <n_k\endSb 
\Big\| \sum_1^k a_i F_{n_i}(\ep) x_{n_i}\Big\| -\ep \\
&\ge \liminf\Sb n_1\to\infty \\ n_1<n_2<\cdots < n_k\endSb 
\frac1K \Big\| \sum_1^k a_i F_{n_i} (\ep) x_{n_i}\Big\|_{\ell_1} -\ep 
\ge {K-\ep\over K} -\ep\ ,
\endalign$$ 
where $\|\cdot\|_{\ell_1}$ refers to the $\ell_1$-norm w.r.t.\ the 
coordinates $(s_i)$. 
Since $\lim_n \|y_n\|=1$ (e.g., use 3)~)  we obtain that $(y_n)$ has a 
spreading model 1-equivalent to the unit vector basis of $\ell_1$.\qed 
\enddemo 

Our argument was motivated by \cite{J}. 

$T$ (see e.g., \cite{CS}, \cite{FJ}) is the completion of the linear space 
of finitely supported real valued sequences under the implicit norm 
$$\|x\| = \|x\|_\infty \vee \sup \biggl\{ \tfrac12 \sum_{i=1}^n \|E_ix\| : 
n\in \nat\ ,\ n\le E_1<\cdots < E_n\biggr\}\ .$$

\demo{Proof of Theorem B} 
We may assume that $X$ has a basis $(b_i)$ which is a block basis of $(t_i)$, 
the unit vector basis for $T$. 
Let $(e_i)$ be a normalized block basis of $(b_i)$ where $e_i$ is a 
$(1+\ep_i)-\ell_1^{m_i}$ average for some sequences $m_i\uparrow \infty$ 
and $\ep_i\downarrow 0$. 
Thus $e_i = (\sum_1^{m_i} \omega_j)/\|\sum_1^{m_i} \omega_j\|$ where 
$(\omega_j)_1^{m_i}$ a normalized block basis of $(b_n)$ which is 
$(1+\ep_i)$-equivalent to the unit vector basis of $\ell_1^{m_i}$. 
By passing to a subsequence we may assume that $(e_i)$ has a spreading model 
$(s_i)$. 

Let $x\in\langle t_i\rangle$ be fixed with $x\le t_n$ and let $k\in\nat$, 
$(a_i)_1^k\subseteq\real$. 
Then 
$$\eqalign{&
\varlimsup \Sb n_1\to\infty\\ n_1<\cdots <n_k\endSb 
\sup \biggl\{ \tfrac12 \sum_{j=1}^\ell \Big\|E_j \biggl( x+\sum_1^k a_i e_{n_i}
\biggr) \Big\| : \ell \le E_1 <\cdots <E_\ell\ ,\ \ell\le n\biggr\} \cr
&\qquad \le \|x\| + \tfrac12 \sum_1^k |a_i|\ .\cr}$$ 
This follows from the fact that since $e_{n_i}$ is a $(1+\ep_{n_i})
-\ell_1^{m_{n_i}}$ average 
$$\eqalign{
&\lim_{i\to\infty} \sup \biggl\{ \tfrac12 \sum_{j=1}^\ell 
\|E_j a_i e_{n_i}\| : \ell \le E_1 <\cdots < E_\ell\ ,\ \ell \le n\biggr\}\cr
&\qquad = \tfrac12 |a_i|\quad \text{(see e.g., \cite{OTW}).}\cr}$$ 
It follows that 
$$\Big\| \sum_1^k a_i s_i\Big\| \le \sup_i \biggl( |a_i| + \tfrac12 
\sum_{i+1}^k |a_j|\biggr)\ .$$ 
However since $\lim_i \|e_i\|_\infty =0$ if $k$ is fixed and $\ep>0$, 
then for $i$ sufficiently large we have for some choice of $E_1<\cdots <E_\ell$ 
that $1= \|e_i\| \le \frac12 \sum_1^\ell \|E_j e_i\| +\ep$ where 
$\ell\le m-k$, $m=\min E_1$. 
This yields that for all $k$, $(a_i)_1^k\subseteq \real$,  
$$\Big\| \sum_1^k a_i s_i\Big\| = \max_i \biggl( |a_i| +\tfrac12 
\sum_{i+1}^k |a_j|\biggr)\ .
\tag 4$$ 

All that remains is to show that Proposition 3.1 applies to $(e_i)$ and 
$(s_i)$. 
$(s_i)$ is 2-equivalent to the unit vector basis of $\ell_1$. 
Set 
$$\align 
x_1 &= \tfrac12 s_1 + \tfrac23 s_2\ ,\cr
x_2 &= \left( \tfrac23\right)^2 s_3 + \left( \tfrac23\right)^2 s_4 
+ \tfrac23 s_5\ ,\cr 
x_3 &= \left(\tfrac23\right)^3 s_6 +\left(\tfrac23\right)^3 s_7 
+ \left( \tfrac23\right)^2 s_8 + \tfrac23 s_9\ ,\ \text{ etc.} 
\endalign$$ 
In general $x_n$ has the same distribution as 
$$ \left( \tfrac23\right)^n s_1 + \left(\tfrac23\right)^n s_2 
+ \left( \tfrac23\right)^{n-1} s_3 +\cdots + \tfrac23 s_{n+1}$$ 
and $(x_n)$ is a block basis of $(s_i)$. 
It is easy to check by (4) that $\|x_n\|=1$ and $\lim_n \|x_n\|_{\ell_1} =2$. 
Also $(x_n)$ is Cauchy in distribution since for $n<m$  
$$d(x_n,x_m) \le \left( \tfrac23\right)^n + \sum_{n+1}^m 
\left(\tfrac23\right)^i +\left(\tfrac23\right)^m\ .
\eqno\qed$$
\enddemo 


\remark{Remark 4.2} 
The above argument yields the following. 
Let $(x_i)$ be a normalized basic sequence having spreading model $(e_i)$ 
equivalent to the unit vector basis of $\ell_1$. 
Let 
$$K\equiv \sup \biggl\{ \sum_1^n |a_i| : \Big\| \sum_1^n a_i e_i\Big\| =1
\biggr\}\ .$$ 
Let  $E_\que$ be the completion of $\langle e_q :q\in \que\rangle$ under 
$\|\sum_1^n a_i e_{q_i}\| = \|\sum_1^n a_i e_i\|$ if $q_1<\cdots <q_n$. 
Suppose there exists $x\in E_\que$ with $\|x\|=1$ and 
$\|x\|_{\ell_1}=K$. 
Then there exists a normalized block basis $(y_i)$ of $(x_i)$ having 
spreading model 1-equivalent to the unit vector basis of $\ell_1$. 
However such an $x$ need not exist (consider $\|x\| =\|x\|_{c_0} 
+\|x\|_{\ell_1}$). 
\endremark 

\remark{Remark 4.3} 
If $|\cdot|$ is any equivalent norm on $T$ then for all $\ep>0$ there 
exists a spreading model of a normalized block basis of $(e_i)$ which 
is $1+\ep$-equivalent to the unit vector basis of $\ell_1$. 
In fact one has \cite{BL, p.43} more generally if $(e_i)$ is a basic 
sequence with spreading model equivalent to the unit vector basis of 
$\ell_1$, then for all $\ep>0$ there exists a spreading model of a normalized 
block basis of $(e_i)$ which is $(1+\ep)$-equivalent to the unit vector basis 
of $\ell_1$.
\endremark 

In \cite{OS} it is proved that if $X$ does not contain $\ell_1$ then $X$ 
can be renormed so that if $(s_i)$ is a spreading model of a normalized 
sequence $(x_n)$ then $\|s_1+s_2\|=1$ implies that $(x_n)$ is not weakly 
null. 
Here we give an explicit renorming of $T$ with this property. 

\proclaim{Proposition 4.4} 
$T$ can be given an equivalent norm $\dotnormm$ satisfying that if $(s_n)$ 
is a spreading model of a normalized block basis of $(T,\dotnormm)$ then 
$$\normm s_1+s_2\normm <2\ .$$
\endproclaim 

First we construct an equivalent norm for $T$. 
Fix $0<q<1/2$ and let $\dotNorm$ be defined on $c_{00}$ by the implicit 
equation 
$$\|x\| = \|x\|_\infty \vee \sup \biggl\{ \frac12 \sum_{i=1}^n 
\|E_ix\| + q\max_{j\le n} \| E_jx\| : n\in\nat\text{ and } n\le E_1<
\cdots < E_n\biggr\}\ .
\tag 1$$ 
As in \cite{FJ} $\dotNorm$ extends to a norm on the completion $X$ of 
$(c_{00},\dotNorm)$ which satisfies (1) for all $x\in X$. 
We postpone the proof of our next proposition which says that $X$ is $T$ 
under an equivalent norm for certain $q$. 

\proclaim{Proposition 4.5} 
$\dotNorm$ is an equivalent norm on $T$ if $(\frac12 +q)^2<\frac12$. 
\endproclaim 

Given $n\in\nat$ and $x\in c_{00}$ set 
$$\align
|x|_n &= \max \biggl\{ \frac12 \sum_{i=1}^n \|E_ix\| : n\le E_1<\cdots E_n
\biggr\}\cr 
\intertext{and} 
\|x\|_n &= \max \biggl\{ \frac12 \sum_{i=1}^n \|E_i x\| + q\max_{j\le n} 
\|E_jx\|: n\le E_1 <\cdots < E_n\biggr\}\ .
\endalign$$ 
As usual we let $(e_i)$ be the unit vector basis of $c_{00}$. 

\proclaim{Lemma 4.6} 
Let $(y_i)$ be a block basis of $(e_i)$ satisfying: 
for all $\ep>0$ there exists $m\in\nat$ with $1\ge \lim_{n\to\infty} |y_n|_m
>1-\ep$. 
Then  $\limsup_{n\to\infty} \|y_n\|>1$. 
\endproclaim 

\demo{Proof} 
By passing to a subsequence we may assume that for all $m\in \nat$, 
$$\lim_{n\to\infty} |y_n|_m = : 1-\ep_m \ \text{ where }\ \ep_m\downarrow0\ .$$

We may also assume that for any $n>m$ there exists sets $m\le E_1^{(m,n)} 
<\cdots < E_m^{(m,n)}$ so that 
$$1-2\ep_m \le\tfrac12 \sum_{j=1}^m \|E_i^{(m,n)} y_n\| \le |y_n|_m \ .
\tag 1$$ 
Furthermore we may assume that there exists a subsequence $M$ of $\nat$ 
such that if $m<m'$ are integers in $M$ then for all $n\ge m'$ and 
$i\le m'$ there exists $j\le m$ with $E_i^{(m',n)} \subseteq E_j^{(m,n)}$. 
Indeed  we can first choose sets to satisfy (1) with lower estimate 
``$1-\frac32 \ep_m$.'' 
If $m_1<\cdots < m_k$ in $M$ have been chosen we then choose $m_{k+1}$ so 
large that the chosen sets for $m_{k+1}$ can be split, if necessary, up to 
$m_1+\cdots + m_k$ additional 
times by first deleting, if necessary, the smallest $m_1+\cdots 
+ m_k$ terms of the form $\|E_\ell^{m_{k+1},n} y_n\|$ at a cost of at most 
$\ep_m/2$. 

Choose $\ep_0>0$ so that $2(\frac12 +q)(1-\ep_0) >1$ and choose $m_0\in M$ 
with $2\ep_{m_0} <\ep_0$. 
Thus for $n>m_0$ 
$$\align
&\left(\tfrac12 +q\right) \sum_{i=1}^{m_0} \|E_i^{(m_0,n)} y_n\| 
> 2\left( \tfrac12 +q\right) (1-\ep_0) \tag 2\cr 
&\qquad > 1 = \lim \Sb m\to\infty \\ m\in M\endSb \lim_{n\to\infty} 
\tfrac12 \sum_{i=1}^{m_0} \sum \Sb j\le m\\ E_j^{(m,n)} \subseteq E_i^{(m_0,n)}
\endSb \|E_j^{(m,n)} y_n\|\ .
\endalign$$

From (2) we see that for some subsequence $M'$ of $M$ we have for  some 
$i_0 \le m_0$ 
$$\left(\tfrac12 +q\right) \lim\Sb n\in M'\\ n\to\infty\endSb 
\|E_{i_0}^{(m_0,n)} y_n\| 
> \lim\Sb m\in M'\\ m\to\infty \endSb 
\lim \Sb n\in M'\\ n\to\infty \endSb \tfrac12 
\sum \Sb j\le m\\ E_j^{(m,n)} \subseteq E_{i_0}^{(m_0,n)} \endSb 
\|E_j^{(m,n)} y_n\|\ .$$ 
Thus 
$$\eqalignno{ 
\limsup_{n\to\infty} \|y_n\| 
& \ge \limsup_{m\to\infty} \limsup_{n\to\infty} 
\biggl[ \left(\tfrac12 +q\right) \|E_{i_0}^{(m_0,n)} y_n\| \cr
&\qquad + \tfrac12 \sum\Sb i\ne i_0\\ i\le m_0\endSb\ \  
\sum \Sb j\le m\\ E_i^{(m,n)} \subseteq E_i^{(m_0,n)}\endSb 
\|E_i^{(m,n)} y_n\| \cr 
& > \limsup_{m\to\infty} \limsup_{n\to\infty} \tfrac12 
\sum_{i=1}^m \|E_i^{(m,n)} y_n\| = 1\ .&\qed\cr}$$
\enddemo 

\demo{Proof of Proposition 4.4} 

Let $\dotNorm$ be the norm given in Proposition 4.5 (for some fixed $q$). 
Let $\dotnormm = \dotNorm +\dotnorm$ where $|x| = \sup \{ \frac12 
\sum_{i=1}^n \|E_ix\| : n\le E_1 <\cdots < E_n\}$, i.e., $|x|=\sup_n |x|_n$. 
$\dotnormm$ is an equivalent norm on $T$. 
Let $(x_n)\subseteq T$ be a $\dotnormm$ semi-normalized block basis of 
$(e_i)$ with a spreading model $(s_n)$. 
We shall prove that $\normm s_1 +s_2\normm <2\normm s_1\normm$. 

If not then we may assume that $\|x_n\|=1$ for all $n$, $\lim_n |x_n|=A$,  
$\lim_{m\to\infty} \lim_{n\to\infty} |x_m+x_n| =2A$ and 
$\lim_{m\to\infty} \lim_{n\to\infty} \| x_m+x_n\|=2$. 
By passing to a subsequence using Ramsey's theorem we may assume  that for all 
$m\in\nat$ there exist $\ell_1(m) \in\nat$ so that 
$\|x_m+x_n\| = \|x_m+x_n\|_{\ell_1(m)}$ whenever $m<n$. 

Also we note that $A<1$. 
Otherwise Lemma~4.6 could be applied to a suitable subsequence of $(x_n)$ 
to yield $\limsup_{n\to\infty} \|x_n\|>1$. 

For each $m<n$ choose $i(m,n) \le \ell_1(m)$ and $\ell_1(m) \le E_1^{(m,n)} 
<\cdots < E_{\ell_1(m)}^{(m,n)}$ with 
$$\|x_m + x_n\| = \left( \tfrac12 +q\right) \|E_{i(m,n)}^{(m,n)} (x_m +x_n)\| 
+ \tfrac12 \sum_{i\ne i(m,n)} \|E_i^{(m,n)} (x_m+x_n)\|$$ 
\enddemo 

By passing again to a subsequence we may assume that we have one of the 
following three cases. 

\remark{Case 1} 
For all $m<n$, $E_{i(m,n)}^{(m,n)} \cap \supp (x_n)=\emptyset$.
\endremark

In this case 
$$\lim_{m\to\infty} \lim_{n\to\infty} \|x_m+x_n\| 
\le \lim_{m\to\infty} \lim_{n\to\infty} (|x_n| + \|x_m\|) = A+1<2\ ,$$ 
a contradiction. 

Case 2 is handled similarly. 

\remark{Case 2} 
For all $m<n$, $E_{i(m,n)}^{(m,n)} \cap \supp (x_m)=\emptyset$.
\endremark 

\remark{Case 3} 
For all $m<n$ 
$$E_{i(m,n)}^{(m,n)} \cap \supp (x_n) \ne \emptyset \ \text{ and }\ 
E_{i(m,n)}^{(m,n)} \cap \supp (x_m)\ne \emptyset\ .$$ 
\endremark 

Since $A<1$ we may assume by similar considerations to those in cases 1) 
and 2) that 
there exists $\ep_0 >0$ with for all $m<n$ 
$$\left( \tfrac12 +q\right) \|E_{i(m,n)}^{(m,n)} x_m\| > \ep_0\ \text{ and }
\ \left( \tfrac12 +q\right) \|E_{i(m,n)}^{(m,n)} x_n\| > \ep_0\ .
\tag 1$$ 
Now for $m<n$, 
$$\align 
\|x_m+x_n\| &\le \tfrac12 \sum_{i=1}^{i(m,n)-1} \|E_i^{(m,n)} x_m\| \tag 2\cr 
&\qquad +\left(\tfrac12 +q\right) \|E_{i(m,n)}^{(m,n)} x_m\| 
+ \left(\tfrac12 +q\right) \|E_{i(m,n)}^{(m,n)} x_n\| \cr 
&\qquad + \tfrac12 \sum_{i=i(m,n)+1}^{\ell_1(m)} \|E_i^{(m,n)} x_n\|\le 2\ .
\endalign$$ 
Since $\lim_{m\to\infty} \lim_{n\to\infty} \|x_m+x_n\|=2$ we may assume 
$$\lim_{m\to\infty}\lim_{n\to\infty}\| E_{i(m,n)}^{(m,n)} (x_m+x_n)\| 
= \lim_{m\to\infty}\lim_{n\to\infty} \left[ \|E_{i(m,n)}^{(m,n)} x_m\| 
+ \|E_{i(m,n)}^{(m,n)} x_n\| \right] \ .$$ 

Thus by once more passing to a subsequence we may assume that for all $m$ 
there exists $\ell_2(m)$ so that if $m<n$ then 
$$\|E_{i(m,n)}^{(m,n)}x_n\| < \|E_{i(m,n)}^{(m,n)} x_n\|_{\ell_2(m)} 
+ {\ep_0\over2} \left( 1-q-\tfrac12\right)\ .$$
Now for each $m$ if $n$ is sufficiently large (e.g., if $\min\supp (x_n) 
\ge \ell_1 (m) + \ell_2(m)$) then 
$$\align 
\|x_n\| & \ge \|E_{i(m,n)}^{(m,n)} x_n\|_{\ell_2(m)} 
+ \tfrac12 \sum_{i=i(m,n)+1}^{\ell_1(m)} \|E_i^{(m,n)} x_n\| \\ 
& \ge \left( 1-q-\tfrac12\right) \|E_{i(m,n)}^{(m,n)} x_n\| 
- {\ep_0\over2} \left(1-q-\tfrac12\right) \\
&\qquad + \left(\tfrac12 +q\right) \|E_{i(m,n)}^{(m,n)} x_n\| 
+ \tfrac12 \sum_{i=i(m,n)+1}^{\ell_1(m)} \|E_i^{(m,n)} x_n\|\ .
\endalign$$ 
As $m\to\infty$, $m<n$ the last line converges to 1 by virtue  of (2). 
From (1) we have $\|E_{i(m,n)}^{(m,n)}x_n\|>\ep_0$. 
Thus the entire expression is, in the limit, 
$$\ge \left( 1-q-\tfrac12\right) {\ep_0\over2} +1 >1\ ,$$ 
a contradiction.\qed

\demo{Proof of Proposition 4.5} 

For $0<c<1$ define the following (implicit) norms on $c_{00}$. 
$$\align
N_c^{(1)} (x) &= \|x\|_\infty \vee \max \biggl\{ \tfrac12 
\sum\Sb i=1\\ i\ne i_0\endSb ^n N_c^{(1)} (E_ix) + cN_c^{(1)} (E_{i_0} x) :\\  
&\hskip1.25truein n\le E_1 <\cdots < E_n\ ,\ 1\le i_0\le n\biggr\}\ .
\endalign$$
Thus the norm $\dotNorm$ in (1) equals $N_c^{(1)}(x)$ if $c=\frac12 +q$. 
$$\align
N_c^{(2)} (x) & = \|x\|_\infty \vee \max \Biggl\{ \tfrac12 
\sum \Sb i=1\\ i\ne i_0\endSb ^n N_c^{(2)} (E_i x) \\
&\qquad + c\max \biggl\{ \sum_{j=1}^m N_c^{(2)} (F_jx) : m\le F_1<\cdots 
< F_m\ ,\\
&\qquad\qquad E_{i_0-1} < F_1\ ,\ F_m <E_{i_0+1}\biggr\} : 
n \le E_1 <\cdots < E_n\ ,\ 1\le i_0 \le n\Biggr\}
\endalign$$
We also define 
$$\align
N^{(2)} (x) & = N_{1/2}^{(2)} (x)\ ,\\
N^{(3)}(x) & = \|x\|_\infty \vee \max \biggl\{ \tfrac12 \sum_{i=1}^n 
N^{(2)} (E_ix) : n\le E_1 <\cdots < E_n\biggr\}\ ,\\
N^{(4)} (x) & = \|x\|_\infty \vee \max \biggl\{ \tfrac12 \sum_{i=1}^{3n} 
N^{(4)} (E_ix) : n\le E_1 < \cdots < E_{3n}\biggr\}\ .
\endalign$$ 

We shall prove that $N_c^{(1)}$ and $N^{(4)}$ are equivalent norms when 
$c^2\le \frac12$ and $c\ge \frac12$. 
This will complete the proof. 
Indeed $(e_i)$, the unit vector basis for $T$ is equivalent to $(e_{3i})$ 
(\cite{CJT} or \cite{CS, p.35}) and for all $(a_i)\in c_{00}$, 
$$N^{(4)} (\sum a_i e_i) = \|\sum a_i e_{3i}\|_T\ .$$ 
\enddemo

We establish three claims. 

\proclaim{Claim 1} 
$N_c^{(1)} (x) \le 2N_{c^2}^{(2)} (x)$ if $c\ge \frac12$.
\endproclaim 

The proof is by induction on $k=\#\,\supp (x)$. 
If $k=1$ the claim is clear. 
Assume that $N_c^{(1)} (x) \le 2N_{c^2}^{(2)} (x)$ if $\#\,\supp (x)\le k$. 
Let $\#\,\supp x=k+1$. 
We may assume that 
$$N_c^{(1)} (x) = \tfrac12 \sum \Sb i=1\\ i\ne i_0\endSb ^n N_c^{(1)} 
(E_ix) + cN_c^{(1)} (E_{i_0}x)$$ 
for some $n\le E_1 <\cdots < E_n$ and $1\le i_0\le n$ with $E_{i_0}x\ne0$ 
and $E_jx\ne 0$ for some $j\ne i_0$. 
Thus we can apply the induction hypothesis to every $E_ix$. 

\remark{Case 1} 
$N_c^{(1)} (E_{i_0}x) = \|E_{i_0}x\|_\infty$.

Then by the induction hypothesis, 
$$\align 
N_c^{(1)} (x) &\le 2\Biggl[ \tfrac12 \sum\Sb i=1\\ i\ne i_0\endSb ^n 
N_{c^2}^{(2)} (E_ix) + \frac c2 \| E_{i_0} x\|_\infty\biggr] \\ 
&\le 2\biggl[ \tfrac12 \sum\Sb i=1\\ i\ne i_0\endSb ^n N_{c^2}^{(2)} 
(E_ix)
+ c^2 \max \biggl\{ \sum_{j=1}^m N_{c^2}^{(2)} (F_jx) : m\le F_1 < 
\cdots < F_m\ ,\\
&\hskip2.0truein E_{i_0} <F_1\ ,\ F_m < E_{i_0+1} \biggr\}\Biggr]\\
&\le 2N_{c^2}^{(2)} (x)\ .
\endalign$$ 
The second inequality uses that $\frac c2 \le c^2$ since $c\ge \frac12$. 
\endremark

\remark{Case 2} 
$$N_c^{(1)} (E_{i_0}x) = \tfrac12 \sum\Sb j=1\\ j\ne j_0\endSb ^m 
N_c^{(1)} (F_jx) + c N_c^{(1)} (F_{j_0} x)$$ 
for some $m\le F_1 <\cdots < F_m$, $j_0\le m$ with $F_{j_0}x\ne0$. 

We may assume (all our norms are 1-unconditional) that $E_{i_0-1} < F_1 < F_m 
< E_{i_0+1}$. 
By the induction hypothesis, 
$$\align  N_c^{(1)}(x) 
& \le 2\biggl[ \tfrac12 \sum \Sb i=1\\ i\ne i_0\endSb ^n N_{c^2}^{(2)} 
(E_ix) + \frac c2 \sum\Sb j=1\\ j\ne j_0\endSb ^m N_{c^2}^{(2)} (F_jx) 
+ c^2 N_{c^2}^{(2)} (F_{j_0} x)\biggr] \\ 
&\le 2\biggl[ \tfrac12 \sum \Sb i=1\\ i\ne i_0\endSb ^n N_{c^2}^{(2)} 
(E_ix) + c^2 \sum_{j=1}^m N_{c^2}^{(2)} (F_jx)\biggr] \\ 
&\le 2N_{c^2}^{(2)} (x)\ ,\ \text{where we have again used that } \frac c2 
\le c^2\ .
\endalign$$ 
This completes the proof of Claim 1.
\endremark

Note that if $\frac12 \le c<1$ and $c^2 \le\frac12$ then 
$N_c^{(1)} (x) \le 2N_{c^2}^{(2)} (x) \le 2N^{(2)}(x)$. 

\proclaim{Claim 2} 
$N^{(3)} (\cdot) \le N^{(4)} (\cdot)$
\endproclaim 

The proof is again by induction on $k=\#\,\supp (x)$. 
If $k=1$ the claim is clear. 
Assume the claim holds for all $x$ with $\#\,\supp (x) \le k$ and 
let $\#\,\supp (x) = k+1$. 
We may suppose that for some $n\le E_1 <\cdots < E_n$, 
$N^{(3)} (x) = \frac12 \sum_{i=1}^n N^{(2)} (E_ix)$. 
For each $i$ either $N^{(2)} (E_ix) = \|E_ix\|_\infty$ or there exist 
integers $s(i) \le n(i)$ and $m(i)$ and sets 
$$n(i) \le E_1^i <\cdots < E_{s(i)-1}^i < F_1^i <\cdots < F_{m(i)}^i 
< E_{s(i)+1}^i <\cdots < E_{n(i)}^i$$ 
with $m(i) \le F_1^i$ so that 
$$N^{(2)} (E_ix) = \tfrac12 \sum_{s=1}^{s(i)-1} N^{(2)} (E_s^i x) + \tfrac12  
\sum_{t=1}^{m(i)} N^{(2)} (F_t^i x) 
+ \tfrac12 \sum_{s=s(i)+1}^{n(i)} N^{(2)} (E_s^i x)\ .$$ 
Thus if we set 
$$E_i^+ = \bigcup_{s> s(i)} E_s^i\ ,\quad 
E_i^- = \bigcup_{s<s(i)} E_s^i\ \text{ and }\ 
E_i^0 = \bigcup_{j=1}^{m(i)} F_j^i$$ 
in the latter case and take $E_i^+ = E_i$ (and $E_i^- = E_i^0 =\emptyset$) 
in the former we have 
$$\align 
N^{(3)} (x) & \le \tfrac12 \sum_{i=1}^n \left( N^{(3)} (E_i^-x) +N^{(3)} 
(E_i^0 x) + N^{(3)} (E_i^+ x)\right) \\ 
&\le \max \biggl\{ \tfrac12 \sum_{i=1}^{3n} N^{(3)} (\tilde E_ix) 
: n\le \tilde E_1 <\cdots < \tilde E_{3n}\biggr\} \\ 
&\le \max \biggl\{ \tfrac12 \sum_{i=1}^{3n} N^{(4)} (\tilde E_i x) : n 
\le \tilde E_1 <\cdots < \tilde E_{3n}\biggr\}\\ 
&\le N^{(4)} (x)\ .
\endalign$$ 
We used the induction hypothesis to obtain the next to last inequality. 

\proclaim{Claim 3} 
$N^{(2)} (\cdot) \le 2N^{(4)} (\cdot)$ 
\endproclaim 

As in Claim 2 we may assume that 
$$\align 
N^{(2)} (x) & = \tfrac12 \sum_{i=1}^{i_0-1} N^{(2)} (E_ix) +\tfrac12 
\sum_{i=i_0+1}^n N^{(2)} (E_ix) +\tfrac12 \sum_{j=1}^m N^{(2)} (F_jx)\\ 
&\le N^{(3)} \biggl( \bigcup \Sb i=1\\ i\ne i_0\endSb ^n E_i x\biggr) 
+ N^{(3)} \biggl( \bigcup_{j=1}^m F_j x\biggr) \\ 
& \le 2N^{(3)} (x)
\endalign$$ 
Thus by Claim 2, the result follows. 

The proposition now follows from our claims: 
$$ \dotNorm_T  \le N_c^{(1)}(\cdot) 
\le 2N_{c^2}{(2)} (\cdot) 
\le 2N^{(2)} (x) 
\le 4N^{(4)} (\cdot)$$ 
where the second inequality uses that 
$c=\frac12 +q$ satisfies $c^2<\frac12$.\qed 

\remark{Remark} 
The proof shows that 
none of our equivalent norms can arbitrarily distort $T$. 
Moreover they all satisfy $|\sum_1^n x_i| \ge \frac12 \sum_1^n |x_i|$ if 
$n\le x_1 <\cdots < x_n$. 
By \cite{OTW, Theorem~6.2} such norms do not arbitrarily distort $T$. 
\endremark

\enddocument